\documentclass{amsart}

\usepackage{amsfonts}
\usepackage{amssymb}
\newcommand{\N}{{\mathbb N}}
\newcommand{\Z}{{\mathbb Z}}

\newcommand{\G }{\Gamma (G, X\cup \mathcal H)}
\newcommand{\GG }{\Gamma (G, X\cup \mathcal H')}

\newcommand{\Hl }{\{ H_\lambda \} _{ \lambda \in \Lambda  }}
\renewcommand{\P}{the property $\rm P_{nai}$}
\newcommand{\PP}{property $\rm P_{nai}$ }

\theoremstyle{plain}
\newtheorem{lem}{Lemma}
\newtheorem{cor}[lem]{Corollary}
\newtheorem{thm}[lem]{Theorem}

\theoremstyle{definition}
\newtheorem*{defn}{Definition}

\begin{document}

\title{Relatively Hyperbolic Groups are $C^*$-simple } %with unique trace}

\author{G. Arzhantseva}
\address{Universit\'{e} de Gen\`{e}ve,
Section de Math\'{e}matiques,
2-4 rue du Li\`{e}vre,
Case postale 64,
1211 Gen\`{e}ve 4, Switzerland}
\email{Goulnara.Arjantseva@math.unige.ch}

\author{A. Minasyan}
\address{Universit\'{e} de Gen\`{e}ve,
Section de Math\'{e}matiques,
2-4 rue du Li\`{e}vre,
Case postale 64,
1211 Gen\`{e}ve 4, Switzerland}
\email{aminasyan@gmail.com}

\begin{abstract} We characterize relatively hyperbolic groups whose
reduced \\{$C^*$-algebra} is simple as those, which have no non-trivial finite
normal subgroups.
%As an easy byproduct we determine relatively hyperbolic groups with the infinite conjugacy classes.
\end{abstract}
\thanks{This work was supported by the Swiss
National Science Foundation Grant $\sharp$~PP002-68627.}
%%2000 MSC:

\maketitle

\section{Introduction}
Let $G$ be a countable discrete group.
We denote by $\ell ^2(G)$ the Hilbert space of square-summable complex-valued functions on $G$
and by $B(\ell ^2(G))$ the algebra of bounded operators on $\ell ^2(G)$.
The group $G$ acts on $\ell ^2(G)$
by means of the left regular representation:
$$
\lambda(g)f(h)=f(g^{-1}h), \ g,h\in G, f\in \ell ^2(G).
$$
The {\it reduced $C^*$-algebra} $C^*_r(G)$ of $G$ is the operator norm closure of the linear span of the
set of operators $\{ \lambda(g) \mid g\in G\}$ in $B(\ell ^2(G))$. It has a unit element and
the canonical trace $\tau : C^*_r(G)\to \mathbb{C}$, given by
$\tau (1)=1$ and $\tau (\lambda(g))=0$ for all $g\in G\setminus\{1\}$.
This $C^*$-algebra reflects analytic properties
of the group $G$.
It plays an important role in non-commutative geometry and, via $K$-theory, in the
Baum-Connes conjecture \cite{val}.

We say that the group $G$ is {\it $C^*$-simple} if its reduced $C^*$-algebra is simple, that is,
it has no non-trivial two-sided ideals.

In 1975 Powers established the $C^*$-simplicity of non-abelian free groups \cite{powers}. Later,
many other examples of $C^*$-simple groups were found. These include
non-trivial free products \cite{pasa},
torsion-free non-elementary Gromov hyperbolic groups \cite{H88,H85}
(more generally, torsion-free non-elementary convergence groups \cite{PH}),
centerless mapping class groups and outer automorphism groups of free groups \cite{bh},
many irreducible Coxeter groups \cite{fen}, etc. A nice overview
of $C^*$-simple groups can be found in \cite{PH}.

The $C^*$-simplicity can be regarded as a strong form of non-amenability:
if $G$ is both amenable and $C^*$-simple then $G$ is reduced to one element~\cite{PH}.
It is a classical result that the existence of free
subgroups in the group implies non-amenability. In \cite{B-C-H}, Bekka, Cowling and de la Harpe
introduced the following ``free-like'' property for finite subsets of a group:
\begin{defn}
A discrete group $G$ is said to have {\it \PP} if
for any finite subset $F$ of $G\setminus \{1\}$ there exists an element $g_0\in G$ of infinite order such that
for each $f \in F$, the subgroup $\langle f,g_0\rangle$ of $G$, generated by $f$ and $g_0$,
is canonically isomorphic to the free product $\langle f\rangle * \langle g_0\rangle$.
\end{defn}
This property guarantees that $C_r^*(G)$ is simple and has a unique normalized
trace~\cite{B-C-H}.
Recall that a {\it normalized trace} on a $C^*$-algebra $A$ with unit is a
linear map $\sigma : A\to \mathbb{C}$ such that $\sigma(1)=1$, $\sigma(a^*a)\ge 0$, and $\sigma(ab)=\sigma(ba)$
for all $a,b$ in $A$.

In the present article our main goal is to characterize $C^*$-simple relatively hyperbolic groups.
The class of relatively hyperbolic groups, that is, groups
hyperbolic with respect to appropriate collections of subgroups, is very large.
It includes Gromov hyperbolic groups and many other examples.
For instance, if $M$ is a complete Riemannian finite--volume manifold of pinched negative sectional curvature,
then $\pi _1(M)$ is hyperbolic with respect to the cusp subgroups
\cite{Bow,Farb}. More generally, if $G$ acts isometrically and
properly   discontinuously on a proper hyperbolic metric space $X$
so that the induced action of $G$ on $\partial X$ is geometrically
finite, then $G$ is hyperbolic relative to  the collection of
maximal parabolic subgroups \cite{Bow}. Groups acting on $CAT(0)$
spaces with isolated flats are hyperbolic relative to the
collection of flat stabilizers \cite{KH}. Algebraic examples of
relatively hyperbolic groups include free products and their small
cancellation quotients \cite{Osin-RHG}, and fully residually free groups
 (or Sela's limit groups) \cite{Dah}.
% , and, more generally, groups acting freely  on $\mathbb R^n$--trees \cite{Gui}.

The notion of a relatively hyperbolic group was originally suggested by Gromov \cite{Gro} and since
then it has been investigated from different points of view
\cite{Bow,Farb,DSO,Thomas}. We use a general approach suggested by Osin in \cite{Osin-RHG}
(see the next section for details) which, when applied to finitely generated groups,
is equivalent to those elaborated by Bowditch~\cite{Bow} and Farb (with the ``bounded coset
penetration" condition)~\cite{Farb}.

Every non-elementary
group $G$ which is hyperbolic relative to a
collection of proper subgroups, or
NPRH group\footnote{This refers to a non-elementary properly relatively hyperbolic group.
A group is {\it elementary} if it has a cyclic subgroup of finite index.} for brevity,
has a {\it maximal finite normal subgroup} denoted by $E_G(G)$~\cite[Lemma 3.3]{SQ}.
The quotient $G/E_G(G)$ is again a NPRH group \cite[Lemma 4.4]{SQ}.
The main result of this paper is the following

\begin{thm} \label{thm:P_nai} Let $G$ be a non-elementary group hyperbolic relative to a collection of
proper subgroups $\Hl$. If $E_G(G)=\{1\}$ then $G$ satisfies {\P}.
\end{thm}

\begin{cor} \label{cor:C-simple} Let $G$ be a NPRH group. Then the following are equivalent.
\begin{itemize}
\item[(i)] The reduced $C^*$-algebra of  $G$
 is simple;
\item[(ii)] The reduced $C^*$-algebra of $G$ has a unique normalized trace;
\item[(iii)] $G$ has infinite conjugacy
classes\footnote{A group $G$ has {\it infinite conjugacy classes}, or, shortly, is icc,
if it is infinite and if all its conjugacy classes distinct from $\{1\}$ are infinite.
A group $G$ is icc if and only if the von Neumann algebra $W^*(\Gamma)$ is  a factor of type $II_1$~\cite[Lemma 5.3.4]{MN}.};
\item[(iv)] $G$ does not have non-trivial finite normal subgroups.
\end{itemize}
\end{cor}

\begin{proof} A discrete $C^*$-simple group can not have a non-trivial amenable normal subgroup \cite[Prop. 2]{Bek-Har00}.
A similar argument shows the same for a discrete group whose the reduced $C^*$-algebra has
a unique normalized trace~\cite[Prop. 2]{Bek-Har00}.
Evidently, a group with infinite conjugacy classes contains no non-trivial finite normal subgroups.
Therefore each of the properties (i) -- (iii) implies (iv).

Suppose, now, that (iv) holds. Then Theorem \ref{thm:P_nai}
and the result of Bekka-Cowling-de la Harpe \cite{B-C-H} mentioned above imply
(i) and (ii).  Property (iii) follows from (i) as noted in ~\cite{PH}.
\end{proof}

Observe that a $C^*$-simple group is icc but the converse is not true in general as there exist amenable
icc groups. However, it is an open problem to find a group $G$ such that $C^*_r(G)$ is simple
with several normalized tracial forms, or a group such that $C^*_r(G)$ has a unique normalized trace
and is not simple~\cite[Prop. 2]{Bek-Har00}.

Every non-elementary Gromov hyperbolic group is a NPRH group with respect to the family consisting of the trivial subgroup.
Therefore Corollary \ref{cor:C-simple} also describes all $C^*$-simple Gromov hyperbolic groups.

Recall that any countable group $G$ has a maximal normal amenable subgroup $R_a(G) \lhd G$ called the {\it amenable radical} of $G$.
Corollary \ref{cor:C-simple} together with the existence of a maximal finite normal subgroup,
imply that for any NPRH group $G$, the quotient $G/E_G(G)$ is $C^*$-simple. Since $R_a(G/E_G(G))=R_a(G)/E_G(G)$ it follows
(see \cite[Prop. 2]{Bek-Har00}) that $R_a(G)=E_G(G)$. We have just obtained

\begin{cor} The amenable radical of a NPRH group $G$ coincides with its maximal finite normal subgroup
$E_G(G)$ and the quotient $G/E_G(G)$ is $C^*$-simple with a unique normalized trace.
\end{cor}
Thus, a NPRH group $G$ is $C^*$-simple if and only if its amenable radical is reduced to one element.
It is worth noticing that it is not yet known whether there exists a countable group $G$
whose amenable radical is trivial but $G$ is {\it not} $C^*$-simple~\cite[Question 4]{PH}.

Easy examples of $C^*$-simple NPRH groups include  non-abelian fully residually free groups
mentioned above~\cite{Dah}. Indeed, since non-abelian free groups satisfy {\P} in an obvious way, a standard
argument\footnote{Consider a finite subset $F \subset G \setminus \{1\}$.
By the assumptions, there exist a non-abelian free group $H$ and a homomorphism $\varphi: G \to H$
such that $\varphi(F) \subset H\setminus\{1\}$. Any subgroup of a free group is free itself, hence we can
suppose that $\varphi$ is surjective. As $H$ satisfies {\P}, we can choose $h \in H$ so that
for every $f \in F$, the elements $\varphi(f)$ and $h$ freely generate a free subgroup $H'\le H$ of rank $2$.
Choose an arbitrary preimage $g_0$ of $h$ in $G$. For any $f \in F$ the restriction of $\varphi$ to the subgroup
$\langle f,g_0 \rangle$ is an isomorphism with $H'$, because a non-trivial element in its kernel would yield a non-trivial
relation between the images $\varphi(f)$ and $\varphi(g_0)=h$. Thus,
$\langle f,g_0 \rangle \cong \langle f \rangle * \langle g_0 \rangle$, and $G$ satisfies {\P}. }
shows the same for $G$.
This method, however, can not be applied to prove Theorem \ref{thm:P_nai}, because there exist $C^*$-simple NPRH groups which
are not fully residually hyperbolic (more generally, which are not which are not limits of Gromov hyperbolic
groups in the space of marked groups -- see \cite{ChG} for the definitions). As an example one can take
a free product $S*{\mathbb Z}$, where $S$ is an infinite finitely presented simple group.

\medskip

{\it Acknowledgments.} The authors would like to thank Pierre de la Harpe
for fruitful discussions.

\section{Relatively hyperbolic groups and their properties}

Let $G$ be a group, $\Hl $ a
fixed collection of subgroups of $G$ (called {\it peripheral
subgroups}), $X$ a subset of $G$. We say that $X$ is a {\it relative
generating set of $G$} with respect to $\Hl $ if $G$ is generated by
$X$ together with the union of all $H_\lambda $. In this situation the group $G$ can be
considered as a quotient of the free product
%\begin{equation}
$$F=\left( \ast _{\lambda\in \Lambda } H_\lambda  \right) \ast F(X),$$
%\label{F} \end{equation}
where $F(X)$ is the free group with the basis $X$.
Let $\mathcal R $ be a subset of $F$ such that the kernel
of the natural  epimorphism $F\to G$ is the normal closure of $\mathcal R $ in the group $F$;
we say that $G$ has {\it relative presentation}
\begin{equation}\label{eq:G}
\langle X,\; \{H_\lambda\}_{\lambda\in \Lambda} \mid R=1,\, R\in\mathcal R\rangle .
\end{equation}
If the sets $X$ and $\mathcal R$ are finite, the relative
presentation (\ref{eq:G}) is said to be {\it finite}.

Define $\mathcal H=\bigsqcup_{\lambda\in \Lambda} (H_\lambda\setminus \{ 1\} )$.
A finite relative presentation \eqref{eq:G} is said to satisfy a {\it linear
relative isoperimetric inequality} if there exists $C>0$ such that, for
every word $w$ in the alphabet $X\cup \mathcal{H}$ (for convenience,  we assume further on that
$X^{-1}=X$) representing the identity in the group
$G$, one has
%\begin{equation} \label{prod}
$$w=_F\prod\limits_{i=1}^k f_i^{-1}R_i^{\pm 1}f_i,$$
%\end{equation}
with the equality in the group $F$, where $R_i\in \mathcal{R}$,
$f_i\in F $, for $i=1, \ldots , k$, and  $k\le C\| w\| $,
where $\| w\|$ is the length of the word $w$.

The group $G$ is called {\it  relatively hyperbolic with respect to a
collection of peripheral subgroups} $\Hl $, if $G$ admits a finite relative
presentation (\ref{eq:G})  satisfying a  linear
relative isoperimetric inequality. This definition is
independent of the choice of the finite generating set $X$ and
the finite set $\mathcal R$ in (\ref{eq:G}) (see \cite{Osin-RHG}).

For a  combinatorial path $p$ in the Cayley graph $\G$ of
$G$ with respect to $X\cup \mathcal H$, we denote by $p_-$, $p_+$ and $\|p\|$
the {\it initial point}, the {\it end point} and the {\it length} correspondingly.
We will write $elem(p)$ for the {\it element} of $G$ represented by the label of $p$.
Further, if $\Omega$ is a subset of $G$ and $g$ belongs to the subgroup $\langle \Omega \rangle \le G$
generated by $\Omega$, then
$|g|_\Omega$  will denote the {\it length of a shortest word} in $\Omega^{\pm 1}$, representing
$g$.
%The corresponding word metric on $G$ is denoted by $\dxh $.

Suppose $q$ is a path in $\G $. Using the terminology from \cite{Osin-RHG}, a subpath $p$ is
called an {\it $H_\lambda $--component} (or, simply, a {\it component}) of $q$,
if the label of $p$ is a word in the alphabet
$H_\lambda\setminus \{ 1\} $ for some $\lambda \in \Lambda $, and $p$ is not contained in a longer
subpath of $q$ with this property.

Two components $p_1, p_2$ of a path $q$ in $\G $ are called {\it
connected} if they are $H_\lambda $--components for the same
$\lambda \in \Lambda $ and there exists a path $c$ in $\G $
connecting a vertex of $p_1$ to a vertex of $p_2$ whose label
entirely consists of letters from $ H_\lambda $. In
algebraic terms, this means that all vertices of $p_1$ and $p_2$
belong to the same coset $gH_\lambda $ for a certain $g\in G$.
 We can always assume $c$ to have length at most $1$, as
every non-trivial element of $H_\lambda $ is included in the set of
generators of $G$.  An $H_\lambda $--component $p$ of a path $q$ is
called {\it isolated } if no other $H_\lambda $--component of
$q$ is connected to $p$.  %A path $q$ is said to be {\it without
%backtracking} if all its components are isolated.

The next statement is often useful in the study of relatively hyperbolic groups.

\begin{lem}[\cite{Osin-RHG}, Lemma 2.27] \label{Omega}
Suppose that a group $G$ is hyperbolic relatively to a collection of
subgroups $\Hl $. Then there exist a finite subset $\Omega
\subseteq G$ and a constant $K\ge 0$ such that the following
condition holds. Let $q$ be a cycle in $\G $, $p_1, \ldots , p_k$
a set of isolated $H_\lambda $--components of $q$ for some
$\lambda \in \Lambda $,  $g_1, \ldots , g_k$ elements of $G$
represented by labels of $p_1, \ldots , p_k$
respectively. Then $g_1, \ldots , g_k$ belong to the subgroup
$\langle \Omega\rangle \le G$  and the word lengths of $g_i$'s
with respect to $\Omega $ satisfy the inequality
$$ \sum\limits_{i=1}^k |g_i|_\Omega \le K\|q\|.$$
\end{lem}

An element $g\in G$ is called {\it hyperbolic} if it is not conjugate to an
element of some $H_\lambda $, $\lambda\in \Lambda $.
%A group $H$ is called {\it elementary} if it has a cyclic subgroup of finite index.
The following
description of elementary subgroups in a relatively hyperbolic group was obtained by Osin.

\begin{lem}[\cite{Osin-ESBG}, Thm. 4.3, Cor. 1.7]\label{Eg}
Let $G$ be a group hyperbolic relatively to a collection of
subgroups $\Hl $.
Let $g$ be a hyperbolic element of infinite order of $G$. Then the
following conditions hold.
\begin{enumerate}
\item The element $g$ is contained in a unique maximal elementary
subgroup $E_G(g)$ of $G$; moreover,
%\begin{equation} \label{eq:elem}
$$E_G(g)=\{ f\in G\; :\;
f^{-1}g^nf=g^{\pm n}\; {\rm for \; some\; } n\in \mathbb N\}.$$
%\end{equation}

\item The group $G$ is hyperbolic relative to the collection
$\Hl\cup \{ E_G(g)\} $.
\end{enumerate}
\end{lem}

A significant restriction on the choice of peripheral subgroups is described in the lemma below.

\begin{lem}[\cite{Osin-RHG}, Thm. 1.4] \label{malnorm}
Suppose that a group $G$ is hyperbolic relative to a collection of
subgroups $\Hl $. Then
\begin{enumerate}
\item[(a)]  For any $g\in G$ and any $\lambda , \mu \in \Lambda $,
$\lambda \ne \mu $, the intersection $H_\lambda^g\cap H_\mu $ is
finite.
\item[(b)] For any $\lambda \in \Lambda $ and any $g\notin
H_\lambda $, the intersection $H_\lambda^g \cap H_\lambda $ is
finite.
\end{enumerate}
\end{lem}

\section{Proof of the main result}
Throughout this section we assume that $G$ is a non-elementary group hyperbolic relatively to a family of
{\it proper} subgroups $\Hl$.

\begin{lem} \label{lem:pow-free-prod} Let $g \in G$ be a hyperbolic element of infinite order satisfying $E_G(g)=\langle g \rangle$.
Then there exists $N\in \N$ such that for any $n \ge N$ and $\lambda \in \Lambda$, the subgroup
$\langle H_\lambda, g^n \rangle \le G$ is canonically isomorphic to the free product $H_\lambda * \langle g^n \rangle$.
\end{lem}

\begin{proof} By Lemma \ref{Eg}, $G$ is hyperbolic relatively to the collection $\Hl \cup \langle g \rangle$.
Lemma~\ref{Omega}
provides a finite subset $\Omega$ of $G$ and a constant $K > 0$ corresponding to this new family of peripheral subgroups.
Since the order of $g$ is infinite, there exists $N \in \N$ such that
\begin{equation} \label{eq:m} g^m \notin \{y \in \langle \Omega \rangle ~:~|y|_\Omega \le 8K\} ~\mbox{ provided } |m| \ge N.
\end{equation}

Fix $n\ge N$ and $\lambda \in \Lambda$. Suppose, by the contrary, that the subgroup of $G$, generated by
$H_\lambda$ and $\langle g^n \rangle$, is not canonically isomorphic to their free product. Then there exist $k \in \N$,
$l_1,\dots,l_k \in n\Z \setminus \{0\}$ and $x_1,\dots,x_k \in H_\lambda \setminus \{1\}$ such that
$$x_1 g^{l_1} x_2 g^{l_2} \cdots x_k g^{l_k}=_G 1.$$
 Define $\mathcal{H}'= \bigsqcup_{\lambda \in \Lambda} (H_\lambda \setminus \{1\}) \sqcup (\langle g \rangle \setminus \{1\})$.
Consider the cycle $q=p_1p_2\dots p_{2k}$ in the Cayley graph $\GG$, where $elem(p_1)=x_1$, $elem(p_2)=g^{l_1}$,
$\dots$, $elem(p_{2k-1})=x_k$, $elem(p_{2k})=g^{l_k}$ and  $\|p_i\|=1$, $i=1,\dots,2k$.

Suppose, at first, that there are two $\langle g \rangle$-components $p_s$ and $p_t$ of $q$, $2 \le s < t \le 2k$,
which are connected.
Then there exists a path $u$ between $(p_t)_-$ and $(p_s)_+$, labelled by an element from $\langle g \rangle$
(in particular, $\|u\|=1$). Without loss of generality, we may assume that $(t-s)$ is minimal.

Consider the cycle $o=p_{s+1}\dots p_{t-1} u$ in $\GG$.
If $t=s+2$, then $elem(p_{s+1})=elem(u^{-1}) \in H_\lambda \cap \langle g \rangle =\{1\}$ (by Lemma \ref{malnorm}),
which would contradict to the choice of $q$. Hence $t \ge s+4$ and $p_{s+2},p_{s+4},\dots,p_{t-2}$ are isolated
$\langle g \rangle$-components of $o$. Applying Lemma \ref{Omega} to the cycle $o$, we achieve
$elem(p_{s+2}), \dots, elem(p_{t-2}) \in \langle \Omega \rangle$ and
$$|elem(p_{s+2})|_\Omega + \dots + |elem(p_{t-2})|_\Omega \le K\|o\|=K(t-s).$$
Since $elem(p_{s+2}), \dots, elem(p_{t-2}) \in \langle g^n \rangle \setminus \{1\}$, we can use the formula \eqref{eq:m}
to obtain $$8K \frac{t-s-2}{2} \le K(t-s), ~\mbox{ consequently }~4 \le \frac{t-s}{t-s-2} \le 3.$$
A contradiction.

Therefore all $\langle g \rangle$-components of $q$ have to be isolated.
Applying Lemma \ref{Omega} to the cycle $q$ leads to a contradiction by the
same argument as before. Thus, the statement is proved.
\end{proof}

Each NPRH group has a hyperbolic element of infinite order \cite[Cor. 4.5]{Osin-ESBG}, therefore we can
use the following simplification of Lemma 3.8 from \cite{SQ}:

\begin{lem} \label{lem:good-elem} Every NPRH group $G$ contains a hyperbolic element $g$ of infinite order
 such that $E_G(g)=\langle g \rangle \times E_G(G)$.
\end{lem}

\begin{proof}[Proof of Theorem \ref{thm:P_nai}]
Consider an arbitrary finite subset $F=\{f_1,\dots,f_m\}$ of $G$. For every $i$, if $f_i \in F$ is an element of finite order,
then we can include the finite subgroup $\langle f_i \rangle$ in the collection of peripheral subgroups $\Hl$, preserving
the relative hyperbolicity of $G$ (for instance, by the characterization of all ``hyperbolically embedded subgroups''
obtained in \cite[Thm. 1.5]{Osin-ESBG}). If $f \in F$ is a hyperbolic element of infinite order, we use Lemma \ref{Eg} to include
the elementary subgroup $E_G(f)$ in $\Hl$.
Thus, further on, we can assume that for each $i \in \{1,\dots,m\}$ there exist $\lambda_i \in \Lambda$ and
$h_i \in G$ such that
\begin{equation} \label{eq:h_i} \langle f_i \rangle \le h_iH_{\lambda_i}h_i^{-1}. \end{equation}

Since $E_G(G)=\{1\}$ by the assumptions, we can apply Lemma \ref{lem:good-elem} to find a hyperbolic element $g \in G$ of
infinite order such that $E_G(g)=\langle g \rangle$. Therefore the elements $g_i=h_i^{-1} g h_i$ satisfy the
assumptions of Lemma \ref{lem:pow-free-prod}, thus there exist $N_i \in \N$ such that the subgroup
$\langle H_{\lambda_i}, g_i^{n_i} \rangle \le G$ is canonically isomorphic to the free product of
$H_{\lambda_i}$ and $\langle g_i^{n_i} \rangle$ for each $n_i \ge N_i$ and $i=1,\dots,m$.
Set $n=\max\{N_i~|~i=1,\dots,m\}$. Formula \eqref{eq:h_i} and the definition of $g_i$ imply that
the subgroup $\langle f_i,g^n \rangle$  of the group $G$ is isomorphic to the free product
$\langle f_i \rangle * \langle g^n \rangle$ for every $i=1,\dots,m$. Hence $G$ satisfies {\P}.
\end{proof}

\end{document}